\documentclass[12pt,oneside]{amsart}

\usepackage{setspace}

\usepackage{amssymb}   

\usepackage{amsmath}

\usepackage{mathrsfs}

\usepackage{stmaryrd}

\usepackage{amsthm}
\usepackage{newlfont}
\usepackage{amscd}
\usepackage{mathtools}
\usepackage{bm}
\usepackage{tikz-cd}
\usepackage{graphicx}

\usepackage{hyperref}
\usepackage{framed}

\usepackage[scr=boondox]{mathalpha}
\usepackage{mathtools}

\usepackage{enumitem}

\usepackage[all]{xy}

\usepackage{textcomp}

\usepackage{amsbsy}

\newtheorem{theorem}{Theorem}

\newtheorem{proposition}[theorem]{Proposition}%
\newtheorem{remark}[theorem]{Remark}
\newtheorem{lemma}[theorem]{Lemma}

\theoremstyle{definition}
\newtheorem{definition}[theorem]{Definition}

\newcommand{\Ainf}{\mathbf{A}_{\mathrm{inf}}}
\newcommand{\F}{\mathbb{F}}

\newcommand{\NN}{\mathcal{N}}

\newcommand{\q}{\mathfrak q}
\newcommand{\pp}{\mathfrak p}

\newcommand{\BdRy}{\mathbf{B}_{\mathrm{dR},y}}

\newcommand{\RN}[1]{%
  \textup{\uppercase\expandafter{\romannumeral#1}}%
}

\makeatletter
\newcommand*{\da@rightarrow}{\mathchar"0\hexnumber@\symAMSa 4B }
\newcommand*{\da@leftarrow}{\mathchar"0\hexnumber@\symAMSa 4C }
\newcommand*{\xdashrightarrow}[2][]{%
  \mathrel{%
    \mathpalette{\da@xarrow{#1}{#2}{}\da@rightarrow{\,}{}}{}%
  }%
}
\newcommand{\xdashleftarrow}[2][]{%
  \mathrel{%
    \mathpalette{\da@xarrow{#1}{#2}\da@leftarrow{}{}{\,}}{}%
  }%
}
\newcommand*{\da@xarrow}[7]{%
  \sbox0{$\ifx#7\scriptstyle\scriptscriptstyle\else\scriptstyle\fi#5#1#6\m@th$}%
  \sbox2{$\ifx#7\scriptstyle\scriptscriptstyle\else\scriptstyle\fi#5#2#6\m@th$}%
  \sbox4{$#7\dabar@\m@th$}%
  \dimen@=\wd0 %
  \ifdim\wd2 >\dimen@
    \dimen@=\wd2 %
  \fi
  \count@=2 %
  \def\da@bars{\dabar@\dabar@}%
  \@whiledim\count@\wd4<\dimen@\do{%
    \advance\count@\@ne
    \expandafter\def\expandafter\da@bars\expandafter{%
      \da@bars
      \dabar@ 
    }%
  }%
  \mathrel{#3}%
  \mathrel{%
    \mathop{\da@bars}\limits
    \ifx\\#1\\%
    \else
      _{\copy0}%
    \fi
    \ifx\\#2\\%
    \else
      ^{\copy2}%
    \fi
  }%
  \mathrel{#4}%
}
\makeatother

\begin{document}
\author{Heng Du}
\address[Heng Du]{Yau Mathematical Sciences Center, Tsinghua University, Beijing 100084, China}
\email{hengdu@mail.tsinghua.edu.cn}

\title[$\Ainf$ has uncountable Krull dimension]{$\Ainf$ has uncountable Krull dimension}

\begin{abstract}
Let $\mathcal{O}_E$ be a complete discrete valuation ring and $R$ be a perfect ring in characteristic $p$, we also assume $R$ is a complete valuation ring whose valuation group is of rank one and non-discrete, we prove the Krull dimension of the ring $W_{\mathcal{O}_E}(R)$ of $\mathcal{O}_E$-Witt vectors over $R$ is at least the cardinality of the continuum.
\end{abstract}
\maketitle

\section{Introduction} 
Let $R$ be a perfect ring over $\F_p$, also assume $R$ is equipped with a rank one and non-discrete valuation. In this paper, we fix $\mathcal{O}_E$ to be a complete discrete valuation ring, define $\Ainf$ to be the ring of $\mathcal{O}_E$-Witt vectors $W_{\mathcal{O}_E}(R)$ over $R$ as in \cite[Sect. 1.2]{FF}, i.e., elements in $\Ainf$ can be regarded as holomorphic functions in variable $\pi$, for a fixed uniformizer $\pi$ of $\mathcal{O}_E$. The main result of this paper is the following:

\begin{theorem}\label{thm:main}
Krull dimension of $\Ainf$ is at least the cardinality of the continuum.
\end{theorem}
Note that the result of the above theorem for equal characteristic case is due to Kang-Park (cf. \cite[Theorem 10]{KP}) that they prove the Krull dimension of $R[\![X]\!]$ is at least the cardinality of the continuum for any rank one non-discrete valuation ring $R$. Their work improved a result of Arnold that shows the Krull dimension of $R[\![X]\!]$ is infinite for such $R$, cf. \cite{Arn}. In the mixed characteristic case, $\Ainf$ was studied in the work of Fontaine and is the core in $p$-adic Hodge theory. The ring structure of $\Ainf$ for general perfectoid rings is studied by Scholze\cite{Sch}, Kedalaya--Liu\cite{KL}, Bhatt--Morrow--Scholze\cite{BMS}, etc, but there are still lots of unknowns, for example, mentioned in \cite{Ked}. Elements in $\Ainf$ should be considered as formal power series in variable ``$\pi$", so $\Ainf$ should share similar properties as $R[\![X]\!]$. There is a conjecture in this direction on the Krull dimension of $\Ainf$ in \cite{Ked} and \cite[Warning 2.24]{Bha}, then proved to be infinite by Lang--Ludwig in \cite{LL} using a similar argument as Arnold. 

To generalize the method of \cite{KP} to the mixed characteristic setting, we rely on a key insight of Fargues and Fontaine: elements of \( \Ainf \) can be interpreted as holomorphic functions on the space \( \lvert Y \rvert \), which parametrizes untilts of \( \operatorname{Frac}(R) \) over \( E \). In particular, \cite{FF} introduces a framework to define the order of vanishing of an element \( f \in \Ainf \) at any point in \( \lvert Y \rvert \). Furthermore, the ``$p$-adic Jensen's formula" established in \cite{FF} provides a powerful tool for estimating this vanishing order by relating it with the multiplicity of slopes in the Newton polygon associated with \( f \). This paper adapts these foundational ideas to construct and analyze elements and ideals in \( \Ainf \).

We begin by reviewing some basic facts of \( \Ainf \) and the space \( \lvert Y \rvert \) in Section~\ref{sec:val}. Section~\ref{sec:vanishingorder} introduces the concept of the vanishing order of elements in \( \Ainf \) at points in \( \lvert Y \rvert \) and uses it to construct a chain of ideals in $\Ainf$. In Section~\ref{sec: zero and NP}, we revisit the framework of the ``p-adic Jensen formula" in a specific case, which serves as a key tool for the proof of Theorem~\ref{thm:main}, presented in Section~\ref{sec:proof of main thm}. Notably, the proof of Theorem~\ref{thm:main} is constructive, and it will rely on the following lemma from commutative algebra, which will be applied in the argument.

\begin{lemma}[{\cite[Theorem 3]{KP}}]\label{lem:main2}
Let $R$ be an integral domain and assume that there is a chain of ideals $\{I_r\}_{r\in (0,1)}$, such that
\begin{itemize}
    \item[(a)] for $0<r<s<1$, we have $I_{s} \subseteq I_{r}$;
    \item[(b)] for each $0<r<1$, there is $g_{r} \in I_{r}$, such that for all $0<r<s$, and all minimal prime $\q$ over $I_{s}$, $g_{r} \notin \q$. 
\end{itemize}
Then, the Krull dimension of $R$ is at least the cardinality of the continuum.
\end{lemma}

\begin{remark}
Lemma~\ref{lem:main2} applies to general \( R \); in particular, \( R \) does not need to be local. The method developed in this paper is expected to be generalized to prove that perfectoid Tate algebras have uncountable Krull dimension, as discussed in \cite{Gar}.
\end{remark}

\section{The ring \texorpdfstring{$\Ainf$}{Ainf} and the space \texorpdfstring{$\lvert Y \rvert$}{|Y|}}\label{sec:val}

Fix a perfect non-discrete valuation ring \( R \) of characteristic \( p \), and let \( v \) be the valuation map \(v\colon R \to \mathbb{R} \cup \{\infty\} \). Fix also a complete discrete valuation field \( E \) with ring of integers \( \mathcal{O}_E \) and residue field \( \mathbb{F}_q \), where \( q = p^h \) for some positive integer \( h \). Let \( \Ainf = W_{\mathcal{O}_E}(R) \) denote the ring of \( \mathcal{O}_E \)-Witt vectors over \( R \). For any uniformizer \( \pi \in \mathcal{O}_E \), one can show that the projection \( \Ainf \to R = \Ainf / (\pi) \) admits a unique multiplicative section \( [-]: R \to \Ainf \), which is independent of the choice of \( \pi \).

Moreover, using the theory of strict \( \pi \)-rings, one can show that every element \( f \in \Ainf \) has a unique \( \pi \)-expansion. That is, any \( f \in \Ainf \) can be uniquely expressed as 
\[
f = \sum_{i \geq 0} [a_i] \pi^i,
\]
where \( a_i \in R \); see \cite[Sect. 1.2]{FF}.

The ring \( \Ainf \) is equipped with a family of Gauss norms. For any non-negative real number \( s \), and for any element \( f = \sum_{i \geq 0}[a_i]\pi^i \in \Ainf \), define
\[
v_s(f) = \inf_{i \geq 0}\{v(a_i) + is\}.
\]

By \cite[Sect. 1.4]{FF}, for \( f \in \Ainf \) and \( s \geq 0 \), the following properties hold:
\begin{itemize}
    \item \( v_s(f) = \infty \) if and only if \( f = 0 \),
    \item \( v_s(fg) = v_s(f) + v_s(g) \),
    \item \( v_s(f+g) \geq \min\{ v_s(f), v_s(g) \} \).
\end{itemize}

Define \( B^{b,+} := \Ainf[\frac{1}{\pi}] \), and let $B^{b} \coloneqq B^{b,+}[\frac{1}{[a]}]$
for any \( a \neq 0 \) in \( R \). The valuations \( v_s \) on \( \Ainf \) extend to \( B^{b} \). For any \( s > 0 \), the function \( \lvert \cdot \rvert_\rho = q^{-v_s(\cdot)} \) defines a norm on \( B^{b} \). Following \cite[Defn. 1.6.2]{FF}, Write \( \rho = q^{-s} \) and let \( B_\rho \) denote the completion of \( B^{b} \) with respect to the norm \( \lvert \cdot \rvert_\rho \).

\begin{definition}[{\cite[Defn. 2.2.1]{FF}}]
An element \( f = \sum_{i \geq 0} [a_i] \pi^i \in \Ainf \) is called \emph{primitive} if \( a_0 \neq 0 \) and there exists some \( i \geq 0 \) such that \( a_i \) is a unit in \( R \). Define \( \operatorname{deg}(f) \), the \emph{degree} of a primitive element \( f \), to be the smallest \( i \) such that \( a_i \) is a unit in \( R \). A primitive element of positive degree is called \emph{irreducible} if it is not a product of primitive elements of positive degree.
\end{definition}

It is straightforward to check that the set of primitive elements of degree \( 0 \) equals \( \Ainf^\times \), the unit group of \( \Ainf \), and it acts on the set of irreducible primitive elements in \( \Ainf \) by multiplication.

\begin{definition}[{\cite[Defn. 2.3.1]{FF}}]
Define \( \lvert Y \rvert \) to be the set of irreducible primitive elements in \( \Ainf \) modulo the action of \( \Ainf^\times \) described above. Let \( \lvert Y \rvert^{\operatorname{deg}=1} \) be the subset of equivalence classes of degree \( 1 \) primitive elements.
\end{definition}

\section{Vanishing order at elements inside \texorpdfstring{$\lvert Y \rvert$}{|Y|}}\label{sec:vanishingorder}

For \( y \in \lvert Y \rvert \), let \( \pp_y \) be the ideal in \( \Ainf \) generated by any element in \( y \). First, recall the following facts.

\begin{theorem}\label{thm:BdR}
For any $y \in \lvert Y \rvert$, let $L_y \coloneqq (\Ainf/\pp_y)[\frac{1}{\pi}]$. Then $L_y$ is a perfectoid field over $E$. Let $\BdRy^+$ denote the $\pp_y[\frac{1}{\pi}]$-adic completion of $\Ainf[\frac{1}{\pi}]$. Then $\BdRy^+$ is a discrete valuation ring with residue field $L_y$ with normalized valuation map given by  $\operatorname{ord}_y: \BdRy^+ \to \mathbb{N} \cup \{\infty\}$.
\end{theorem}
\begin{proof}
The first part is from (2) of \cite[Thm. 3.3.1]{FF}, and the second part is from Sect. 3.4.1 of \textit{loc. cit.}
\end{proof}

\begin{definition}
For any $f \in \Ainf$, define the vanishing order $\operatorname{ord}_y(f) \in \mathbb{N} \cup \{\infty\}$ of $f$ at $y \in \lvert Y \rvert$ to be the image of $f$ under the composition of the natural map $\Ainf \to \BdRy^+$ with $\operatorname{ord}_y$ in Theorem~\ref{thm:BdR}.
\end{definition}

Fix a sequence of elements $\{y_n\}$ in $\lvert Y \rvert$. For any real number $r$, let $\lfloor r \rfloor$ denote the greatest integer less than or equal to $r$.

\begin{definition}\label{defn:chain of ideals}
For any $r > 0$, let $I_r$ be the subset of $\Ainf$ defined by
\[
I_r \coloneqq \{ f \in \Ainf \mid \exists \text{ real number } c > 0 \text{ such that } 
\operatorname{ord}_{y_n}(f) \geq c n^{\lfloor rn \rfloor} \text{ for } n \gg 0\}.
\]
\end{definition}

Using Theorem~\ref{thm:BdR}, it is immediate to show that the following holds.

\begin{lemma}
The family $\{I_r\}_{r \in (0,1)}$ defines a chain of ideals in $\Ainf$ satisfying condition $(a)$ in Lemma~\ref{lem:main2}.
\end{lemma}

\section{Zeros and Newton polygons}\label{sec: zero and NP}
For \( f = \sum_{n \gg -\infty} [a_n] \pi^n \in B^b \), let \( \NN(f) \) denote the Newton polygon of \( f \). Recall that \( \NN(f) \) is defined as the convex, piecewise-linear function from \( \mathbb{R} \) to \( \mathbb{R} \cup \{\infty\} \), determined by the boundary of the convex hull of the set \( \{i, v(a_i)\}_{i \gg -\infty} \). For a convex, piecewise-linear function \( \mathcal{F} \) from \( \mathbb{R} \) to \( \mathbb{R} \cup \{\infty\} \) that is not identically equal to \( \infty \), \( x \) is called a \emph{node} of \( \mathcal{F} \) if \( \mathcal{F}(x) < \infty \) and \( \mathcal{F} \) is not differentiable at \( x \). 

Note that all nodes of \( \NN(f) \) are integers. The \emph{negatives} of the slopes of the affine segments of \( \NN(f) \) are referred to as the \emph{slopes} of \( \NN(f) \). For example, under this convention, $\NN(f)$ is a decreasing convex function when $f\in \Ainf$, and the the \emph{slopes} of \( \NN(f) \) are \emph{nonnegative}.

For each \( I = [q^{-r_1}, q^{-r_2}] \subset (0,1) \), let \( \NN_I(f) \) denote the Newton polygon obtained by ``removing the affine segments of \( \NN(f) \) with slopes not contained in \( [-r_1, -r_2] \)". We refer the reader to \cite[Definition 1.6.18]{FF} for the precise definition of \( \NN_I(f) \), where \( \NN_I(f) \) is denoted by \( \mathscr{Newt}^0_I(f) \).

In the above definition, note that the case \( r_1 = r_2 = s \) is allowed, and in this case, \( \NN_I(f) \) will be denoted by \( \NN_{\rho}(f) \) with \( \rho = q^{-s} \). Fix $\rho\in (0,1)$ in the rest of this section, we now summarize some key properties of \( \NN_{\rho}(f) \) and the ring \( B_{\rho} \) that will be used in Section~\ref{sec:proof of main thm}.

\begin{proposition}\label{prop about zero and NP}
\begin{enumerate}
    \item Let \( \{f_n\}_{n\geq 1} \) be a sequence in \( B^b \) that converges to \( f \in B_{\rho} \). Then \( \NN_{\rho}(f_n) = \NN_{\rho}(f_{n+1}) \) for \( n \gg 0 \).
    
    \item An element \( f \in B^b \) is invertible in \( B_{\rho} \) if the slopes of \( \NN(f) \) are not equal to \( -\log_q \rho \).
    
    \item When \( r > 0 \) is contained in the valuation group \( v(\operatorname{Frac}(R)^\times) \subset \mathbb{R} \), \( B_{q^{-r}} \) is a principal ideal domain with all maximal ideals are of the form \( \pp_y B_{q^{-r}} \) for some \( y \in \lvert Y \rvert \). Moreover, for any \( f \neq 0 \in \Ainf \), \( f \) admits a factorization \( f = u \xi_1^{m_1} \cdots \xi_l^{m_l} \) in \( B_{q^{-r}} \), with \( u \in B_{q^{-r}}^\times \), \( \{\xi_i\} \) generators of distinct maximal ideals \( \pp_{y_i} B_{q^{-r}} \) in \( B_{q^{-r}} \), and $m_i$ positive integers.
    
    \item The factorization in part (3) is unique up to multiplications by units in $B_{q^{-r}}$. Moreover, for \( y_i \), \( m_i \), and \( f \) as in part (3), \( \operatorname{ord}_{y_i}(f) = m_i \).
\end{enumerate}
\end{proposition}

\begin{proof}
Part (1) and part (2) follow from \cite[Cor. 1.6.19 and Prop. 1.6.25]{FF}. The first half of part (3) is a restatement of \cite[(2) in Thm. 3.5.1]{FF}. The second half of part (3) and part (4) can be deduced from \cite[Thm. 3.5.5]{FF}, using the definition of \( \operatorname{div}(f) \) given on page 147 (after Prop. 3.5.10) and part (2) of Thm. 3.5.1 in \textit{loc. cit.} 
\end{proof}

Part (1) of the above proposition implies that if a convergent sequence \( \{f_n\}_{n \geq 1} \) in \( \Ainf \) converges to some element \( f \in \Ainf \), then the slopes of \( \NN(f) \) and their multiplicities are equal to those of \( \NN(f_n) \) for sufficiently large \( n \).

\section{Proof of the main theorem}\label{sec:proof of main thm}

In this section, we assume that \( E \) is of characteristic \( 0 \), and we will construct, for each \( r \in (0,1) \), an element \( g_r \) in \( \Ainf \) such that \( \operatorname{ord}_{y_n}(g_r) = n^{\lfloor rn \rfloor} \) for a sequence \( \{y_n\} \) in \( \lvert Y \rvert \). If it can be further shown that \( \{I_r\}_{r \in (0,1)} \), as defined in Definition~\ref{defn:chain of ideals}, along with \( \{g_r\} \), satisfies the conditions in Lemma~\ref{lem:main2}, then Theorem~\ref{thm:main} will follow in the case when \( E \) is of characteristic \( 0 \). But as noted previously, the equal characteristic case has already been addressed in \cite{KP}.

\begin{proof}[{Proof of Theorem~\ref{thm:main} assuming Lemma~\ref{lem: gr}.}] Following the above discussion, it is enough to show that \(\{I_r\}_{r \in (0,1)}\), together with \(\{g_r\}\), satisfies the conditions in Lemma~\ref{lem:main2}, and the argument is the same as in \cite[Lemma 9]{KP}, which we briefly recall here. Assume the contrary: there exist \( r, s \in (0,1) \) with \( s > r \), such that \( g_r \) is contained in a minimal prime ideal \(\q\) over \( I_s \). Localizing at \(\q\) and modulo \( I_s \), \( g_r \) becomes nilpotent. This implies that \( g_r^l h \in I_s \) for some \( h \in \Ainf \setminus \q \) and \( l \geq 1 \). A computation shows that \(\operatorname{ord}_{y_n}(h) \geq cn^{\lfloor sn \rfloor} - ln^{\lfloor rn \rfloor}\) for $n\gg 0$ and some $c>0$, in particular \(\operatorname{ord}_{y_n}(h) \geq \frac{c}{2}n^{\lfloor sn \rfloor}\) for $n\gg 0$, which contradicts the assumption that \( h \) is \emph{not} in \( I_s \).
\end{proof}

\medskip
\noindent
\textbf{Construction of $y_n$.}
Since the valuation of \( R \) is assumed to be non-discrete, there exists a sequence of elements \( \{b_n\}_{n \in \mathbb{N}} \) in \( R \), such that \( \{v(b_n)\} \) is a \emph{strict} decreasing sequence of positive real numbers satisfying \( v(b_n) \leq n^{-n-1}2^{-n} \) for all \( n \). It is direct to check that
\begin{equation}\label{eq:1}
   \sum_{n \geq 1} n^{\lfloor rn \rfloor + 1} v(b_n) \leq 1 
\end{equation}
for all \( r \in (0,1) \).

Let \( y_n \) be the primitive element of degree \( 1 \) defined by \( \pi - [b_n] \). Since the valuations of \( b_n \) are assumed to be distinct, the \( y_n \) are distinct.

Fix \( \{b_n\}_{n \geq 1} \) as above. Also, fix \( d \in R \) with \( v(d) = 1 \) for the rest of this section.

\medskip
\noindent
\textbf{Construction of $g_r$.} For any \( r \in (0,1) \) and $n\geq 1$, define
\[
g_{r,n} = [d]^2 \prod_{i=1}^n \left(1 - \frac{\pi^i}{[b_i]^i}\right)^{i^{\lfloor ri \rfloor}}.
\]
From Eq.~\eqref{eq:1}, one has \( g_{r,n} \in \Ainf \) for all \( n \). It is straightforward to check that \( g_{r,n+1} \equiv g_{r,n} \mod \pi^{n+1} \) for \( n \geq 1 \). In particular, the sequence \( \{g_{r,n}\} \) defines an element \( g_r \in \Ainf \).

Similarly, for all \( j \geq 1 \), one can define \( \hat{g}^j_r \in \Ainf \) by the sequence
\[
\hat{g}^j_{r,n} = [d] \prod_{i=1, i \neq j}^n \left(1 - \frac{\pi^i}{[b_i]^i}\right)^{i^{\lfloor ri \rfloor}}.
\]
Moreover, it is straightforward to check that \( g_r \) admits a factorization \(g_r =\left([d]\left(1 - \frac{\pi^j}{[b_j]^j}\right)^{j^{\lfloor rj \rfloor}}\right) \hat{g}^j_r \) in \( \Ainf \) for all \( r \in (0,1) \) and \( j \geq 1 \), by verifying a similar equation for \( g_{r,n} \) and \( \hat{g}^j_{r,n} \).

\begin{lemma}\label{lem: gr}
For all \( r \in (0,1) \) and \( n \geq 1 \),
\[
\operatorname{ord}_{y_n}(g_r) = n^{\lfloor rn \rfloor}.
\]
\end{lemma}

\begin{proof}
Let \( s_i \coloneqq v(b_i) \). Using part (1) of Proposition~\ref{prop about zero and NP}, the slopes of \( \NN(g_r) \) are \( \{s_i\}_{i\geq 1} \) for all \( r \in (0,1) \), and the slopes of \( \NN(\hat{g}^j_r) \) are \( \{s_i\}_{i\geq 1,i \neq j} \).

Let \( \rho_n \coloneqq q^{-s_n} \). By part (2) of Proposition~\ref{prop about zero and NP}, \( \hat{g}^n_r \) is a unit in \( B_{\rho_n} \). By parts (3) and (4) of the same proposition, verifying \( \operatorname{ord}_{y_n}(g_r) = n^{\lfloor rn \rfloor} \) reduces to checking that, after writing
\[
[d]\left(1 - \frac{\pi^n}{[b_n]^n}\right)^{n^{\lfloor rn \rfloor}} = w ([b_n] - \pi)^{n^{\lfloor rn \rfloor}},
\]
it holds that \( \operatorname{ord}_{y_n}(w) = 0 \). 

This follows from the fact that
\[
\frac{[b_n]^n - \pi^n}{[b_n] - \pi} \neq 0 \mod (\pi - [b_n]),
\]
when viewed as an element in \( L_{y_n} \). Here, the assumption that \( \operatorname{char}(E) = 0 \) is used.

\end{proof}

\medskip
\noindent
\textbf{Acknowledgments.} We thank Jaclyn Lang and Judith Ludwig for their paper on the related topic. We thank Pavel \v Coupek, Kiran Kedlaya, Tong Liu, Judith Ludwig, Linquan Ma, Dongming She, and Yifu Wang for their feedback on the early versions of the paper. We also thank Kiran Kedlaya for suggesting this more accurate statement of the main theorem. Additionally, we extend our gratitude to the anonymous referee for pointing out a critical mistake in an earlier version of the manuscript, which led to the development of the entirely new proof presented here. This work was supported by the National Key R\&D Program of China (No. 2023YFA1009703) and the Beijing Municipal Natural Science Foundation (Youth Program, No. 1254044).

\bibliographystyle{amsplain}
\bibliography{library}

\end{document}